\documentclass[11pt,reqno]{amsart}
\usepackage[english]{babel}    \usepackage[latin1]{inputenc}   \usepackage[T1]{fontenc}     \usepackage[french]{minitoc}
\usepackage[nice]{nicefrac}    \usepackage{latexsym,amsfonts}  
\usepackage{graphics}    \usepackage{ulem}       \usepackage{hhline}    \usepackage{dsfont}    \usepackage{mathrsfs}
\usepackage{fancyhdr}    \usepackage{amsmath}    \usepackage{amssymb}   \usepackage{rotating}  \usepackage{fancybox}
\usepackage{color}       \usepackage{colortbl}   \usepackage{setspace}  \usepackage{enumerate} \usepackage{amsthm}
\usepackage{multicol}    
\usepackage{amsthm}    \usepackage{varioref}  \usepackage{textcomp}
\usepackage{lmodern}     \usepackage{mathpazo}   \usepackage{euscript}  \usepackage[pdftex]{hyperref}
\newcommand{\C}{\mathbb C}

\begin{document}

\title[]{\textsf{An integral representation for Folland's fundamental solution of the
sub-Laplacian on Heisenberg groups $\Bbb{H}^{n}$} }
\author{Allal GHANMI $\&$ Zouha\"{i}r MOUAYN }
\address{{\bf (A.G.)} \quad  Department of Mathematics,  Faculty of  Sciences, P.O. Box 1014,
   Mohammed V University,  Agdal,  10 000 Rabat - Morocco  }
\address{{\bf (Z.M.)} \quad Department of Mathematics, Faculty of Sciences and Technics (M'Ghila), P.O.
Box 523, Sultan Moulay Slimane University, 23 000 B\'{e}ni Mellal - Morocco}

\date{\today}
\maketitle

\begin{abstract}
We prove that the Folland's fundamental solution for the sub-Laplacian on
Heisenberg groups \cite{Folland73} can also be derived form the resolvent
kernel of this sub-Laplacian \cite{AskMou}. This provides us with a new integral representation for this fundamental solution.
\end{abstract}

\quad

The Heisenberg group $\Bbb{H}^{n}$ is the nilpotent Lie group whose underlying manifold is $\Bbb{C}^{n}\times \Bbb{R}$ with coordinates $(z,\tau)=(
z_{1},\cdots,z_{n};\tau) $ and whose group law is
\begin{equation}
\left( z,\tau \right) .\left( w,s\right) = \left(z+w,\tau +s+2\Im \left\langle z,w\right\rangle\right),  \label{law}
\end{equation}
where $\left\langle z,w\right\rangle =\sum\limits_{j=1}^{n}z_{j}\overline{w_{j}}$. Letting $z_j=x_j+iy_j\in \C,$ then, $x_{1},\cdots,x_{n}$; $y_{1},\cdots,y_{n}$ and $\tau$ are real coordinates on $\Bbb{H}^{n}$. We set
\begin{equation}
\frac{\partial }{\partial z_{j}}=\frac{1}{2}\left( \frac{\partial }{\partial x_{j}}-i\frac{\partial }{\partial y_{j}}\right) , \qquad
\frac{\partial }{\partial\overline{z}_{j}}=\frac{1}{2}\left( \frac{\partial }{\partial x_{j}}+i\frac{\partial }{\partial y_{j}}\right)
.\label{2}
\end{equation}
The operator
\begin{equation}
\Delta _{\Bbb{H}^{n}}
:=\sum\limits_{j=0}^{n}\left[ -\frac{\partial ^{2}}{\partial z_{j}\partial \overline{z}_{j}}-\left| z_{j}\right| ^{2}\frac{
\partial ^{2}}{\partial t^{2}}+i\frac{\partial }{\partial \tau }
\left( z_{j}\frac{\partial }{\partial z_{j}}-\overline{z}_{j}\frac{\partial }{\partial \overline{z}_{j}}\right) \right]
 \label{SubLap}  
\end{equation}
is left-invariant and is subelleptic of order $\frac{1}{2}$ at each point $(z,\tau)$ of $\Bbb{H}^{n}$ (see Theorem 1 in \cite{Folland73}). Also, in \cite{Folland73}, G.B. Folland has used the analogy with the fact that $|x|^{2-m}$ is (a constant multiple of) the fundamental solution of
Laplacian on $\Bbb{R}^{m}$ with source at zero  to prove that the
sub-Laplacian in \eqref{SubLap} admits a fundamental solution with source at zero of the form
\begin{equation}
G_{0}^{F}\left( z,\tau \right) =c_{n}\left( \left| z\right| ^{4}+\tau
^{2}\right) ^{-\frac{n}{2}}; \qquad \left( z,\tau \right) \in \Bbb{H}^{n},  \label{FFS4}
\end{equation}
where
\begin{equation}
 c_{n}=\frac{2^{n}\Gamma ^{2}\left( \frac{n}{2}\right) }{\pi ^{n+1}}. \label{cst}
\end{equation}
In other words
\begin{equation}
\left( \Delta _{\Bbb{H}^{n}}\left[ \phi \right] ,G_{0}^{F}\right) =\phi(0)   \label{6}
\end{equation}
for any function $\phi \in C_{0}^{\infty }\left( \Bbb{H}^{n}\right)$ (see 
 \cite{Folland73,FollandStein74}).

Here, our aim is to derive the Folland's fundamental solution  in \eqref{FFS4} 
from the resolvent kernel of the sub-Laplacian in \eqref{SubLap} on the 
Heisenberg group  $\Bbb{H}^{n}$. This provides us with a new integral representation of $G_{0}^{F}$. Namely, we have the following \\

\noindent {\bf Theorem.} {\it 
The Folland's fundamental solution in \eqref{FFS4} can also be expressed as
\begin{equation}
G_{0}^{F}\left( z,\tau \right) =  
 \frac{2^{n+1} \Gamma\left( \frac{n}{2}\right)}{\pi^{n+1}}
\int_{0}^{+\infty }x^{n-1}e^{-x|z| ^{2}}\Psi \left( \frac{n}{2},n;2x|z| ^{2}\right) \cos \left(
\tau x\right) dx  \label{FFS} 
\end{equation}
in terms of the Tricomi $\Psi$-function.}\\

To prove \eqref{FFS}, we recall first that for $\zeta \in \Bbb{C}$ with $
\Re(\zeta) <0$, the resolvent operator of $\Delta _{\Bbb{H}^{n}}$ has the
form
\begin{equation}
\left( \zeta -\Delta _{\Bbb{H}^{n}}\right) ^{-1}\left[ \varphi \right] \left( z,t\right)
=\int_{\Bbb{H}^{n}}\mathcal{R}\left( \zeta ;\left(z,\tau \right) ,\left( w,s\right) \right) \varphi \left( w,s\right) d\mu(w,s),
\label{ResOp}
\end{equation}
where the kernel function has been obtained in \cite[p.4, Theorem 3.2]{AskMou} as
\begin{align}
\mathcal{R}\left( \zeta ;\left( z,\tau \right) ,\left( w,s\right) \right)
=\frac{-2^n}{\pi ^{n+\frac{1}{2}}}\int_{0}^{+\infty } & x^{n-1}\Gamma\left( \frac{n}{2}-\frac{\zeta }{2x}\right) \Psi \left( \frac{n}{2}-\frac{%
\zeta }{2x},n;2x\left| z-w\right| ^{2}\right)   \label{kernelFct} \\ 
&\times \exp \left( -x\left| z-w\right| ^{2}\right) \cos \left( x\left(
\tau -s\right) +2x\Im \left\langle z,w\right\rangle \right) dx.
\nonumber
\end{align}
The Tricomi $\Psi$-function can be defined as a linear combination of two
$\, {_1\digamma_1}$-sums (\cite[p.56]{Tricomi54}):
\begin{equation}
\Psi \left( a,c;\xi \right) :=\left(\frac{\Gamma(1-c)}{\Gamma(a-c+1)} 
+ \frac{\Gamma(c-1)}{\Gamma(a)}\xi ^{1-c}\right) {_1\digamma_1}(a,c;\xi),   \label{Tricomi5}
\end{equation}
where
\[
\, {_1\digamma_1}\left( a,c;\xi \right) =\frac{\Gamma(c) }{\Gamma(a)} \sum\limits_{j=0}^{+\infty }
\frac{\Gamma(a+j) }{\Gamma(c+j)} \frac{\xi ^{j}}{j!}
\]
is the confluent hypergeometric function defined with the usual condition $c\neq 0,-1,-2, \cdots$.\\
In the limit $\zeta \to 0$ in \eqref{ResOp}, we obtain a
right inverse of the operator in \eqref{SubLap} as
\begin{equation}
\Delta _{\Bbb{H}^{n}}^{-1}\left[ \varphi \right] \left( z,\tau \right)
=\int_{\Bbb{H}^{n}}\mathcal{R}\left( 0;\left( z,\tau \right) ,\left(
w,s\right) \right) \varphi \left( w,s\right) d\mu(w,s) ,  \label{Inv6}
\end{equation}
$d\mu$ being the Lebesgue measure on $\Bbb{H}^{n}$, or equivalently
\begin{equation}
\mathcal{R}_{0}:=-\mathcal{R}\left( 0;(z,\tau),(w,s) \right)
\end{equation}
 is a Green kernel for $\Delta _{\Bbb{H}^{n}}$ as pointed out in \cite[p.6, Remark 3.3]{AskMou}.
 Now, to establish a connection between the integral kernel $\mathcal{R}_{0}$ and the
Folland's fundamental solution, we proceed by computing the integral
\begin{align}
\mathcal{R}_{0}
 =\frac{2^n}{\pi^{n+\frac{1}{2}}}\int_{0}^{+\infty}x^{n-1} 
 & \Gamma\left( \frac{n}{2}\right) \Psi \left( \frac{n}{2},n;2x|z-w| ^{2}\right) \label{Int5} \\ 
 &\times \exp \left( -x\left| z-w\right| ^{2}\right) \cos \left( x\left(
\tau -s\right) +2x\Im \left\langle z,w\right\rangle \right) dx. \nonumber
\end{align}
For this, let us rewrite \eqref{Int5} as
\begin{equation}
\mathcal{R}_{0}=
\frac{2^n}{\pi^{n+\frac{1}{2}}}\int_{0}^{+\infty }x^{n-1}\left[\Gamma\left( \frac{n}{2}\right) \Psi \left( \frac{n}{2},n;\mu x\right)\right]
  e^{-\frac{\mu }{2}x}\cos \left( x\theta \right) dx , \label{Int6}
\end{equation}
where 
we have set
\begin{equation}
\mu :=2|z-w|^{2} \qquad \mbox{and} \qquad \theta =\left( \tau -s\right) +2\Im \left\langle z,w\right\rangle . \label{mutheta}
\end{equation}
Following Tricomi \cite[p.90]{Tricomi54}, there is an integral representation for the $\Psi$-function in \eqref{Tricomi5} of the form:
\begin{equation}
\Gamma\left( a\right) \Psi \left( a,c;u\right) =\int_{0}^{+\infty
}e^{-at}\exp \left( -\frac{u}{e^{t}-1}\right) \left( 1-e^{-t}\right) ^{-c}dt.
\label{7}
\end{equation}
For the parameters $a=\frac{n}{2}$, $c=n$ and $u=\mu x$, the integral in \eqref{Int6} takes the form
\begin{equation}
\mathcal{R}_{0}=
\frac{2^n}{\pi^{n+\frac{1}{2}}}\int_{0}^{+\infty }x^{n-1}\left[
\int_{0}^{+\infty }e^{-\frac{n}{2}t}\exp \left( -\frac{\mu x}{e^{t}-1}%
\right) \left( 1-e^{-t}\right) ^{-n}dt\right] e^{-\frac{\mu }{2}x}\cos
\left( x\theta \right) dx . \label{Int8}
\end{equation}
Intertwining the integrals, we rewrite \eqref{Int8} as
\begin{equation}
\mathcal{R}_{0}=
\frac{2^n}{\pi^{n+\frac{1}{2}}}\int_{0}^{+\infty }\frac{e^{-\frac{n}{2}t}%
}{\left( 1-e^{-t}\right) ^{n}}\left[ \int_{0}^{+\infty }x^{n-1}\exp \left( -%
\frac{\mu x}{2}\coth \left(\frac{t}{2}\right)\right) \cos \left( x\theta \right)
dx\right] dt.  \label{Int9}
\end{equation}
To calculate the integral
\begin{equation}
 \int_{0}^{+\infty }x^{n-1}\exp \left( -%
\frac{\mu x}{2}\coth \left(\frac{t}{2}\right)\right) \cos \left( x\theta \right)
dx , \label{Int91}
\end{equation}
 we make appeal to the identity (\cite[p.498]{GR})
\begin{equation}
\int_{0}^{+\infty }x^{\nu -1}e^{-\alpha x}\cos (\theta x)dx=
\Gamma\left(\nu \right) \left( \alpha ^{2}+\theta ^{2}\right) ^{-\frac{\nu }{2}}\cos
\left( \nu \arctan\left( \frac{\theta }{\alpha }\right) \right) \label{10}
\end{equation}
when $\Re(\nu) >0$, $\Re(\alpha) > |\Im(\theta)|$ are fulfilled. In our case $\nu =n$ and $\alpha =\frac{\mu }{2}\coth \left(\frac{t}{2}\right)$, and therefore \eqref{Int9} reads
\begin{align}
\mathcal{R}_{0}=
\frac{2^n\Gamma(n)}{\pi^{n+\frac{1}{2}}}  \int_{0}^{+\infty }%
\frac{e^{-\frac{n}{2}t}}{\left( 1-e^{-t}\right) ^{n}}&\left( \left( \frac{\mu}{2}
 \coth \left(\frac{t}{2}\right)\right) ^{2}+\theta ^{2}\right) ^{-\frac{n}{2}} \label{Int11}\\
 & \times \cos\left( n \arctan\left( \frac{\theta }{\frac{\mu }{2}\coth \left(\frac{t}{2}\right)}%
\right) \right) dt.  \nonumber
\end{align}
After some calculations, we see that equation \eqref{Int11} takes the form
\begin{align}
\mathcal{R}_{0}=
\frac{\Gamma(n)}{\pi^{n+\frac{1}{2}}}\int_{0}^{+\infty }
&\left( \sinh \left(\frac{t}{2}\right)\right)^{-n}
\left( \frac{\mu}{2} \coth\left(\frac{t}{2}\right)\right) ^{-n} \label{Int12} \\
& \times \left( 1+\left( \frac{\theta }{\frac{\mu }{2}\coth
\left(\frac{t}{2}\right)}\right) ^{2}\right) ^{-\frac{n}{2}}\cos \left( n\arctan\left(
\frac{\theta }{\frac{\mu }{2}\coth \left(\frac{t}{2}\right)}\right) \right) dt. \nonumber 
\end{align}
Making the change variable $\rho =\frac{\theta }{\frac{\mu }{2}\coth \frac{t%
}{2}},$ then \eqref{Int12} becomes
\begin{equation}
\mathcal{R}_{0}=
\frac{2^n\Gamma(n)}{\mu^{n-1} \theta\pi^{n+\frac{1}{2}}} \int_{0}^{\frac{2\theta }{\mu } }
\left( 1-\left( \frac{\mu }{2\theta }\right) ^{2}\rho ^{2}\right) ^{\frac{n}{2}-1}
\left( 1+\rho ^{2}\right) ^{-\frac{n}{2}}\cos \left( n\arctan(\rho )\right) d\rho .  \label{Int13}
\end{equation}
A second change of variable, $\arctan\left( \rho\right) =u$, shows that \eqref{Int13} can be reduced to
\begin{equation}
\mathcal{R}_{0}=
\frac{2^n\Gamma(n)}{\mu^{n-1} \theta\pi^{n+\frac{1}{2}}}
\int_{0}^{\arctan(\frac{2\theta }{\mu })}
\left( \cos^{2}(u)-\left( \frac{\mu }{2\theta }\right) ^{2}\sin ^{2}u\right) ^{\frac{n}{2}-1}\cos (nu)du . \label{Int14}
\end{equation}
Next, by a trigonometrical linearization, we can rewrite  \eqref{Int14} as
\begin{align}
\mathcal{R}_{0}=
\frac{2^{\frac{n}{2}+1}\Gamma(n)}{\mu^{n-1} \theta\pi^{n+\frac{1}{2}}}
&\left( 1+\left( \frac{\mu }{2\theta }\right) ^{2}\right) ^{\frac{n}{2}-1}  \label{Int15} \\ 
&\times \int_{0}^{\arctan(\frac{2\theta }{\mu })}\left( \cos 2u-\frac{%
1-\left( \frac{2\theta }{\mu }\right) ^{2}}{1+\left( \frac{2\theta }{\mu }%
\right) ^{2}}\right) ^{\frac{n}{2}-1}\cos (nu)du .\nonumber
\end{align}
Setting $\beta =\left( \frac{2\theta }{\mu }\right) ^{2}$ and making the change
of variable $2u=\kappa $ in \eqref{Int15} give that
\begin{equation}
\mathcal{R}_{0}=
\frac{2^{\frac{n}{2}}\Gamma(n)}{\mu^{n-1} \theta\pi^{n+\frac{1}{2}}}
\left( \frac{\beta+1}{\beta }\right) ^{\frac{n}{2}-1}\int_{0}^{2\arctan(%
\frac{2\theta }{\mu })}\left( \cos (\kappa) -\frac{1-\beta }{1+\beta }\right)
^{\frac{n}{2}-1}\cos \left(\frac{n}{2}\kappa \right)d\kappa .  \label{Int16}
\end{equation}
Now, by putting
\begin{equation}
\cos (\varepsilon) =\frac{1-\beta }{1+\beta },  \label{Int17}
\end{equation}
 we can check that
\begin{equation}
\varepsilon =\arccos \left(\frac{1-\beta }{1+\beta }\right)=
\arccos\left(\frac{1-\left( \frac{\mu }{2\theta }\right) ^{2}}{1+\left( \frac{\mu }{2\theta }\right) ^{2}}\right)=2\arctan\left( \frac{2\theta }{\mu }\right).   \label{Int18}
\end{equation}
We are now in position to apply the identity (\cite[p.406]{GR}):
\begin{equation}
\int_{0}^{\varepsilon }\left( \cos (x)-\cos (\varepsilon) \right) ^{\nu -%
\frac{1}{2}}\cos (ax)dx=\sqrt{\frac{\pi }{2}}\left( \sin (\varepsilon) \right)
^{\nu }\Gamma\left( \nu +\frac{1}{2}\right) P_{a-\frac{1}{2}}^{-\nu }\left(
\cos (\varepsilon) \right)   \label{Int19}
\end{equation}
with the conditions $\Re(\nu) >-\frac{1}{2}$, $a>0$ and $0<\varepsilon <\pi$ and where
\begin{equation}
P_{\lambda }^{\nu }(x):= \frac1{\Gamma(1-\nu)}\left(\frac{x+1}{x-1} \right)^{\nu/2} {_2\digamma_1}\left(-\lambda,\lambda +1;1-\nu; \frac{1-x}2\right) \label{LegendreFct}
\end{equation}
denotes the Legendre function of the first kind \cite[p.959]{GR}. In our setting, the integral in \eqref{Int16} reads
\begin{align}
\int_{0}^{2\arctan(\frac{2\theta }{\mu })}\left( \cos (\kappa) -\frac{1-\beta }{1+\beta }\right) ^{\frac{n-1}{2}-\frac{1}{2}}
& \cos \left(\frac{n}{2}
\kappa \right)d\kappa  \label{Int20} \\
& =\sqrt{\frac{\pi }{2}}\left( \sin (\varepsilon) \right) ^{\frac{n-1}{2}}\Gamma\left( \frac{n}{2}\right) P_{\frac{n-1}{2}}^{-\left(
\frac{n-1}{2}\right) }\left( \cos (\varepsilon) \right) .  \nonumber
\end{align}
Hence, using the Gegenbauer representation (\cite[p.969]{GR}):
\begin{equation}
P_{\sigma }^{-\sigma }\left( \cos (\varepsilon) \right) =\frac{1}{\Gamma\left( 1+\sigma \right) }\left( \frac{1}{2}\sin (\varepsilon) \right) ^{\sigma
}  \label{Int22}
\end{equation}
with $\sigma =\frac{n-1}{2}$, we can write the right hand side in \eqref{Int20} as
\begin{align}
\sqrt{\frac{\pi }{2}}\left( \sin (\varepsilon) \right) ^{\frac{n-1}{2}}\Gamma\left( \frac{n}{2}\right) P_{\frac{n-1}{2}}^{-\left( \frac{n-1}{2}\right)
}\left( \cos (\varepsilon) \right)
&=\left( \frac{1}{2}\right) ^{\frac{n-1}{2}}\sqrt{\frac{\pi }{2}}\frac{\Gamma\left( \frac{n}{2}\right) }{\Gamma\left( \frac{n+1}{2}\right) }\left( \sin
(\varepsilon) \right) ^{n-1}
\nonumber
\\&
=\left( \frac{1}{2}\right) ^{\frac{n}{2}}\sqrt{\pi }\frac{\Gamma\left( \frac{n}{2}\right) }{\Gamma\left( \frac{n+1}{2}\right) }%
\left( 1-\cos ^{2}(\varepsilon) \right) ^{\frac{n-1}{2}} . \label{24}
\end{align}
Returning back to \eqref{Int16}, keeping in mind the expression of $\varepsilon$ given through \eqref{Int18}, we obtain that
\begin{align}
\mathcal{R}_{0}&=\frac{2^{\frac{n}{2}}\Gamma(n)}{\mu^{n-1} \theta\pi^{n+\frac{1}{2}}}
\left( \frac{\beta+1}{\beta }\right) ^{\frac{n}{2}-1}\left(\left( \frac{1}{2}\right) ^{\frac{n}{2}}\sqrt{\pi }\frac{\Gamma\left( \frac{n}{2}\right) }{\Gamma\left( \frac{n+1}{2}\right) }\left(  \frac{4\beta }{(1+\beta)^2 }\right) ^{\frac{n-1}{2}}\right)
  \label{26}
\\&
=\frac{2^{n-1}}{\mu^{n-1} \theta\pi^{n }}
 \frac{\Gamma(n)\Gamma\left( \frac{n}{2}\right) }{\Gamma\left( \frac{n+1}{2}\right) }
\left( \frac{\beta+1}{\beta }\right) ^{-\frac{n}{2}}.
 \end{align}
Now, replacing $\beta$ by $\left( \frac{2\theta }{\mu }\right) ^{2}$ and using the expressions of $\mu$ and $\theta$ as in \eqref{mutheta}, we arrive at
\begin{align}
\mathcal{R}_{0}
&=  \frac{ \Gamma(n)\Gamma\left( \frac{n}{2}\right) }{\pi^n\Gamma\left( \frac{n+1}{2}\right) }
\left( \left(|z-w| ^{2}\right) ^{2}+\left( \left( \tau -s\right) +2\Im %
\left\langle z,w\right\rangle \right) ^{2}\right) ^{-\frac{n}{2}}.  \label{**}\\
&= \frac{ 2^{n-1}\Gamma^2\left( \frac{n}{2}\right)}{\pi^{n+\frac12}}
\left( \left(|z-w| ^{2}\right) ^{2}+\left( \left( \tau -s\right) +2\Im %
\left\langle z,w\right\rangle \right) ^{2}\right) ^{-\frac{n}{2}}.  \label{***}
\end{align}
The last equality follows using Legendre's duplication formula (\cite[p.896]{GR}):
\[
\Gamma\left( \xi \right)
\Gamma\left( \xi +\frac{1}{2}\right)= 2^{1-2\xi}\sqrt{\pi }\Gamma(2\xi)
\]
for $\xi =\frac{n}{2}$. Therefore, we assert that
\begin{align}
\mathcal{R}_{0}  =   \frac{ 2^{n-1}\Gamma^2\left( \frac{n}{2}\right)}{\pi^{n+\frac12}c_n}
 G_{0}^{F}\left( \left( z,\tau \right) \circ \left( w,s\right) ^{-1}\right) ,
 \label{FFSp1}
\end{align}
where the constant $c_n$ is as in \eqref{cst}. In particular, for $\left( w,s\right) =\left( 0,0\right) $, Equation \eqref{FFSp1} 
reduces further to
\begin{equation}
\mathcal{R}_{0}
  =  \frac{ 2^{n-1}\Gamma^2\left( \frac{n}{2}\right)}{\pi^{n+\frac12}c_n} G_{0}^{F}\left(  z,\tau \right)
 = \frac{\sqrt\pi}{2} G_{0}^{F}\left(  z,\tau \right)
 .
 \label{Int*}
 \end{equation}
Finally, by combining \eqref{Int6} and \eqref{Int*}, we get the announced result of the theorem.\\

\end{document}